\documentclass[11pt,a4paper,reqno]{amsart}
\RequirePackage[margin=1in,footskip=30pt]{geometry}
\usepackage{p01preamble}
\usepackage{p01bibformat}
\title{Fourier Multipliers on Quasi-Banach Orlicz Spaces and Orlicz Modulation Spaces}
\author{Albin Petersson }
\date{2026-09-09}
\address{Department of Mathematics,
Linn{\ae}us University, V{\"a}xj{\"o}, Sweden}
\email{albin.petersson@lnu.se}

\begin{document}
\begin{abstract}
    We find that if a Fourier multiplier is continuous from \(L^{\Phi_1}\) to \(L^{\Phi_2}\),
then it is also continuous from \(M^{\Phi_1,\Psi}\) to \(M^{\Phi_2,\Psi}\),
where \(\Phi_1,\Phi_2,\Psi\) are quasi-Young functions and \(\Phi_1\) fulfills
the \(\Delta_2\)-condition.
This result is applied to show that Mihlin's Fourier multiplier theorem and
H\"ormander's improvement hold in certain Orlicz modulation spaces.
Lastly, we show that the Fourier multiplier with symbol \(m(\xi) = e^{i \mu(\xi)}\),
where \(\mu\) is homogeneous of order \(\alpha\), is bounded on
quasi-Banach Orlicz modulation spaces of order \(r\), assuming \(r\in\big(d/(d+2),1\big]\)
and \(\alpha\in\big(d(1-r)/r, 2\big]\).
\end{abstract}
\maketitle

\section{Introduction}\label{sec:intro}
Fourier multipliers are obtained by performing a Fourier transform, multiplying by
a suitable function (which we call the symbol), and lastly performing an
inverse Fourier transform.
With \(T\) as the operator and \(m\) as the symbol, we formally write this as
\[
    T(f) = (\ms F^{-1} \circ m\cdot \circ\ms F)(f).
\]
Since \(\partial_j f(x) = \ms F^{-1} [i\xi \cdot \hat{f}](x)\)
(where \(\hat{f} = \ms F[f]\)), partial differential operators with constant coefficients
are Fourier multipliers, and because of this relationship between \(T\) and its symbol
\(m\), we write \(T=m(D)\).
Evidently, Fourier multipliers naturally appear when solving partial differential equations.
A fundamental question is whether a symbol \(m\) will give rise to a bounded operator \(m(D)\)
on a certain function space.
In this paper, we investigate such boundedness conditions on (quasi-Banach) Orlicz spaces
and Orlicz modulation spaces.


Orlicz spaces, initially introduced by W. Orlicz \autocite{Orlicz1932}, are a
generalization of Lebesgue spaces which are, roughly speaking, obtained by replacing
the integrand \(|f(x)|^p\) in the expression \(\nm{f}[L^p]^p\) with \(\Phi(|f(x)|)\),
where \(\Phi\) is a certain type of convex function.
In general, \(\rho_\Phi(f) = \int \Phi(|f(x)|) dx\) is not a norm (for instance,
it may fail to be homogeneous), hence one has to define the norm in a different way.
There are several (equivalent) ways to do this. Here, we use the Luxemburg norm
\autocite[43]{Luxemburg1955} given by
\[
    \nm{f}[L^\Phi] = \inf\sets*{\lambda>0}{\rho_\Phi(f/\lambda) \le 1}.
\]
Thus, the Orlicz space \(L^\Phi\) is the space of measurable functions \(f\) such that
\(\nm{f}[L^\Phi]\) is finite.
If \(\Phi(t) = t^p\), \(1\le p < \infty\), then
\(\nm{\cdot}[L^\Phi] = \nm{\cdot}[L^p]\) and \(L^\Phi = L^p\).
Other examples of Orlicz spaces include \(L^{\Phi_1}\) and \(L^{\Phi_2}\), where
\(\Phi_1(t) = t \log(1+t)\) and \(\Phi_2(t) = e^{t^2} - 1\).


By replacing \(L^p\) norms with \(L^\Phi\) norms in the definitions of other spaces,
we obtain different Orlicz type spaces.
Among these are the so-called Orlicz modulation spaces, where the usual \(L^{p,q}\) norm
is replaced by the mixed Orlicz norm \(L^{\Phi,\Psi}\). 
The usefulness of such Orlicz modulation spaces is illustrated by the following example.
Let \(E_\fy\) be the entropy functional given by
\[
    E_\fy (f) = - \iint_{\rd[2d]} |V_\fy f(x,\xi)|^2 \log |V_\fy f(x,\xi)|^2 \, dxd\xi
    + \nm{V_\fy f}[L^2]^2 \log \nm{V_\fy f}[L^2]^2.
\]
Here, \(V_\fy\) is the short-time Fourier transform and \(\fy\) is a window function.
This functional appears when dealing with kinetic energy estimates in quantum systems
(cf.~\autocite{LiebSolovej1991}).
While \(E_\fy\) is not continuous on \(M^2 = L^2\), it is continuous on \(M^\Phi\),
where \(\Phi\) is a Young function satisfying
\(\Phi(t) = -t^2 \log t\) for \(0\le t\le e^{-3/2}\).
Moreover, \(M^\Phi\) is a dense subset of \(L^2\), and for any \(p<2\),
\(M^p \subseteq M^\Phi\), so that \(M^\Phi\) is, in some sense, a better setting for
the study of this functional than \(M^p\) is for any \(p<2\).
See \autocite[Section 3]{GumberRanaToftUester2024} for more details.


Here, we consider Fourier multipliers on (quasi-Banach) Orlicz spaces and
Orlicz modulation spaces, fulfilling various conditions.
In particular, we obtain boundedness properties for Fourier multipliers with symbol
\(m(\xi) = e^{i \mu(\xi)}\), where \(\mu\) is homogeneous of order \(\alpha \le 2\).
Such operators appear when solving certain evolution equations.
Consider the initial value problem
\[
    \partial_tf(t,x) = i \mu(D) f(t,x)
\]
where \(f(0,x) = f_0(x)\).
For example, the free Dirac and Schrödinger equations can be described in this way, where
\(\mu\) is homogeneous of order \(1\) and \(2\), respectively.
Formally, the solution to the equation is given by \(f(t,x) = e^{i t \mu} f_0(x)\), hence
concerns of existence of solutions to the equation are directly linked to
the boundedness of the aforementioned operator \(m(D)\).


In \zcref{sec:prelim}, we introduce the necessary notations.
\zcref{sec:fmult} is divided into two parts.
In \zcref{subsec:general}, we use a Marcinkiewicz type interpolation result
(cf.~\autocite[Theorem 5.1]{LiuWang2012}) to extend results about the continuity of
Fourier multipliers on \(L^p\) spaces to \(M^\Phi\) spaces (\zcref{cor:mihlin,cor:horm}).
We also show that Fourier multipliers which are bounded on (quasi-Banach) Orlicz spaces
\(L^\Phi\) are also bounded on (quasi-Banach) Orlicz modulation spaces \(M^{\Phi,\Psi}\)
(\zcref{thm:ezproof}).
In \zcref{subsec:quasibanach}, we focus on quasi-Banach Orlicz modulation spaces and
obtain sufficient conditions for the Fourier multiplier with symbol \(m(\xi) = e^{i\mu(\xi)}\),
where \(\mu\) is homogeneous of a certain order, to be bounded on \(M^{\Phi,\Psi}\) with
\(\Phi\) and \(\Psi\) as quasi-Young functions (\zcref{thm:thm1benyi}).
\section{Preliminaries}\label{sec:prelim}
We write \(d\) for dimension, and we let
\[
    \mb N = \{0,1,2,3,\dots\}
\]
be the set of natural numbers.
For multi-indices \(\alpha, \beta \in \nd\), meaning \(\alpha=(\alpha_1,\dots,\alpha_d)\),
\(\beta=(\beta_1,\dots,\beta_d)\) where \(\alpha_j,\beta_j \in \mb N\), \(j=1,\dots,d\),
we write \(\alpha\le \beta\) to mean \(\alpha_j\le \beta_j\) for every \(j=1,\dots,d\).
Further, we let \(|\alpha| = \alpha_1 + \dots + \alpha_d\), and for \(r\in\rd\), we let
\(\alpha^r = \alpha_1^{r_1} \dots \alpha_d^{r_d}\).
Moreover, for \(j,k\in\zd\) we let
\begin{equation}\label{eq:Qj}
    \begin{aligned}
        Q_{r}(j) & =j+[0,r]^d,        & Q(j) & = Q_1(j), \text{ and } \\
        Q(j,k)   & = Q(j)\times Q(k).
    \end{aligned}
\end{equation}

The Fourier transform of \(f\) is denoted \(\ms F[f]\) or \(\hat{f}\), and is given by
\[
    \ms F[f](\xi) = (2\pi)^{-d/2} \int_{\rd} f(x) e^{-i \scal{x}{\xi}} dx.
\]


Recall that a \emph{quasi-norm} of order \(r\in(0,1]\), or an \emph{\(r\)-norm}, to
the vector space \(\mc B\) is a functional \(\Lambda\) on \(\mc B\) for
which the following holds:
\begin{enumerate}
    \item \(\Lambda(f) \ge 0\) for all \(f\in \mc B\) and
          \(\Lambda(f) = 0\) if and only if \(f=0\);

    \item \(\Lambda(\alpha f) = |\alpha| \rho(f)\)
          for all \(f\in\mc B\), \(\alpha\in \mb C\);

    \item \(\Lambda(f+g)^r \le \Lambda(f)^r + \Lambda(g)^r\) for all \(f,g\in\mc B\).
\end{enumerate}
Although quasi-norms are typically defined in a different way, we use the terms
\qu{quasi-norm of order \(r\)} and \qu{\(r\)-norm} interchangeably, and justify this choice
by the Aoki-Rolewicz theorem (cf.~\autocite{Aoki1942,Rolewicz1957}).


A \emph{quasi-Banach space of order \(r\)} or \emph{\(r\)-Banach space} is a
complete quasi-normed space, meaning, it is complete with respect to
the topology induced by a quasi-norm of order \(r\).
For more information about quasi-Banach spaces, see \autocite{LoristNieraeth2024}.


\subsection{Lebesgue type spaces}


For \(p\in(0,\infty]\), let \(L^p(\rd)= L^p\) denote the usual Lebesgue space with norm
\[
    \nm{f}[L^p(\rd)] \equiv \nm{f}[L^p] \equiv
    \begin{cases}
        \left( \int_{\rd} |f(x)|^p d x \right)^{\frac 1 p}, & p<\infty,
        \\[1ex]
        \displaystyle{\esssup_{x\in \rd}} |f(x)|,           & p=\infty,
    \end{cases}
\]
where \(\fdef\) is a Lebesgue measurable function.
The norm \(\nm{\cdot}[L^p]\) simultaneously imposes decay and growth conditions on
the functions in \(L^p\) which depend on the variable \(p\).
Using different such conditions with respect to different variables, we arrive at
the definition of mixed norm Lebesgue spaces, which we recall below.


\begin{defn}
    The \emph{mixed norm Lebesgue space} \(L^{p,q}(\rd[2d])\) consists of all
    Lebesgue measurable functions \(\fdef[f][\rd[2d]]\) such that
    \[
        \nm{f}[L^{p,q}(\rd[2d])] \equiv\nm{f}[L^{p,q}]
        \equiv
        \nm{f_2^p}[L^q]
    \]
    is finite, where
    \[
        f_2^p(y) = \nm{f(\cdot,y)}[L^p].
    \]
\end{defn}


For a window function \(\fy\in \ms S(\rd)\setminus\{0\}\), we let
\(\fdef[V_\fy][\ms S'(\rd)][\ms S'(\rd[2d])]\) denote the \emph{short-time Fourier transform}
given by
\[
    V_\fy f (x,\xi) =
    (2\pi)^{-d/2}\int f(y) \overline{\fy(y-x)} e^{-i \scal{y}{\xi}} \, d y.
\]

We recall further the following definition.


\begin{defn}
    The \emph{modulation space} \(M^{p,q}(\rd)\) consists of all
    tempered distributions \(f\in\ms S'(\rd)\) such that
    \[
        \nm{f}[M^{p,q}(\rd)] \equiv \nm{f}[M^{p,q}] \equiv \nm{V_\fy f}[L^{p,q}]
    \]
    is finite, where \(\fy\in\ms S(\rd)\setminus\{0\}\) is a window function.
\end{defn}


\begin{rem} \label{rem:equivModNorm}
    Observe that the norm
    \[
        \nm{f}[M^{p,q}]^*\equiv \nm{\tilde{V}_\fy f}[L^{p,q}],
    \]
    with
    \[
        \tilde{V}_\fy f (x,\xi) =
        (2\pi)^{-d/2} \int f(y + x) \overline{\fy(y)} e^{-i \scal{y}{\xi}} d y,
    \]
    is equivalent to \(\nm{f}[M^{p,q}]\), where \(\fy \in \ms S\setminus\{0\}\).
\end{rem}


By letting \(L_*^{p,q}(\rd[2d])\) consist of all
Lebesgue measurable functions \(\fdef[f][\rd[2d]]\) such that
\[
    \nm{f}[L_*^{p,q}(\rd[2d])] \equiv \nm{f}[L_*^{p,q}] \equiv \nm{f_1^q}[L^p]
\]
is finite, where
\[
    f_1^q(x) = \nm{f(x,\cdot)}[L^q],
\]
we obtain the following definition.


\begin{defn}
    The space \(W^{p,q}(\rd)\) consists of all
    tempered distributions \(f\in\ms S'(\rd)\) such that
    \[
        \nm{f}[W^{p,q}(\rd)] \equiv \nm{f}[W^{p,q}] \equiv \nm{V_\fy f}[L_*^{p,q}]
    \]
    is finite, where \(\fy\in\ms S(\rd)\setminus\{0\}\) is a window function.
\end{defn}


Analogous to the \(L^p\) spaces, we recall that
the \emph{discrete Lebesgue spaces} \(\ell^{p}(\zd)\) consist of
sequences \(a=\{a(j)\}_{j\in\zd}\) for which
\[
    \nm{a}[\ell^p(\zd)] \equiv \nm{a}[\ell^p] \equiv
    \begin{cases}
        \left(\sum_{j\in\zd} |a(j)|^p \right)^{\frac{1}{p}}, & 0<p<\infty \\
        \sup_{j\in\zd} |a(j)|,                               & p=\infty
    \end{cases}
\]
is finite. On this topic, we recall further.


\begin{defn}
    The \emph{discrete mixed Lebesgue space} \(\ell^{p,q}(\zd[2d])\) consists of
    sequences \(a=\left\{a(j,k)\right\}_{j,k\in\zd}\) such that
    \[
        \nm{a}[\ell^{p,q}(\zd[2d])] \equiv \nm{a}[\ell^{p,q}] \equiv \nm{a_2^p}[\ell^q]
    \]
    is finite, where \(a_2^p(k) = \nm{a(\cdot,k)}[\ell^p]\), \(k\in\zd\).
\end{defn}


Parallel to \(L_*^{p,q}(\rd[2d])\), we let \(\ell_*^{p,q}(\zd[2d])\) consist of
all sequences \(a=\left\{a(j,k)\right\}_{j,k\in\zd}\) such that
\[
    \nm{a}[\ell_*^{p,q}(\zd[2d])] \equiv \nm{a}[\ell_*^{p,q}] \equiv \nm{a_1^q}[\ell^p]
\]
is finite, where \(a_1^q(j) = \nm{a(j,\cdot)}[\ell^q]\), \(j\in\zd\).

To finish our excursion into definitions of Lebesgue spaces, modulation spaces, and
Wiener amalgam spaces, we have the following.


\begin{defn}
    Let \(\mc B = \ell^{p,q}(\zd[2d])\) or \(\mc B =\ell_*^{p,q}(\zd[2d])\), \(0<p,q\le \infty\).
    The \emph{Wiener space} \(W^r(\mc B)\) consists of
    all functions \(\fdef[F][\rd[2d]]\) such that
    \[
        \nm{F}[W^{r}(\mc B)]\equiv \nm{a}[\mc B] < \infty,
    \]
    where
    \[
        a(j,k) = \nm{F}[L^r(Q(j,k))], \quad j,k\in\zd,
    \]
    with \(Q(j,k)\) as in \eqref{eq:Qj}.
\end{defn}


\begin{rem}\label{rem:modprop}
    Note that, for any window function \(\fy\in \ms S\setminus\{0\}\), there exists
    a constant \(C>0\) such that
    \[
        C^{-1}\nm{V_\fy \cdot}[L^{p,q}]
        \le \nm{V_\fy \cdot}[W^\infty(\ell^{p,q})]
        \le C \nm{V_\fy \cdot}[L^{p,q}]
    \]
    so that the norms are equivalent when restricted to short-time Fourier transforms.
    Note also that this is independent of the choice of window function \(\fy\),
    in the sense that different choices for \(\fy\) lead to equivalent norms
    (cf.~\autocite[Theorem 3.1]{GalperinSamarah2004}).
    In a similar manner, \(\nm{V_\fy \cdot}[L_*^{p,q}]\) and
    \(\nm{V_\fy\cdot}[W^\infty(\ell_*^{p,q})]\) are also equivalent.
    See \autocite{GalperinSamarah2004} (in particular Theorem 3.3) for more details.
\end{rem}

Lastly, we recall the definition of a Fourier multiplier, which plays a pivotal role
in this paper.


\begin{defn}
    Let \(m\in \ms S'(\rd)\). Then the \emph{Fourier multiplier}
    \(\fdef[m(D)][\ms S(\rd)][\ms S'(\rd)]\)  with \emph{symbol} \(m\) is given by
    \[
        m(D) = \ms F^{-1}\circ(m\cdot)\circ\ms F.
    \]
\end{defn}

Observe that if, for example, \(m\in L^p(\rd)\), \(1\le p < \infty\), then
\[
    m(D)f(x) =
    (2\pi)^{-d/2}\int_{\rd}  m (\xi) \hat f (\xi) e^{i \scal{x}{\xi}} d \xi,
    \quad f\in\ms S(\rd).
\]


\subsection{Quasi-Young functions and quasi-Orlicz spaces}


We recall the following definitions of Young functions and quasi-Young functions.


\begin{defn}\label{Def:YoungFunc}
    Let \(\Phi\) be a function from \([0,\infty)\) to \([0,\infty]\).
    Then \(\Phi\) is called a \emph{Young function} if
    \begin{enumerate}
        \item \(\Phi\) is convex,

              \vrum

        \item \(\Phi (0)=0\),

              \vrum

        \item \(\Phi (t)<\infty\) for some \(t>0\),

              \vrum

        \item \(\lim\limits_{t\rightarrow\infty} \Phi (t)=+\infty\).
    \end{enumerate}
\end{defn}


\begin{defn}
    A function \(\Phi\) from \([0,\infty)\) to \([0,\infty]\)
            is called a \emph{quasi-Young function} if
            there exists \(r\in (0,1]\) such that \(t \mapsto \Phi(t^{1/r})\) is a Young function.
    The largest such \(r\) is called the \emph{order} of \(\Phi\).
\end{defn}

Note that a quasi-Young function must be increasing.
The concept of a quasi-Young function is explored in \autocite{ToftUesterNabizadehOeztop2022}, but
can also be found under the name \(s\)-convex \(N\)-function \autocite[43]{RaoRen1991}
(with the additional assumptions that \(\lim_{t\rightarrow 0^+} \frac{\Phi(t)}{t}=0\) and
\(\lim_{t\rightarrow \infty} \frac{\Phi(t)}{t}=\infty\)).

We will briefly consider Lebesgue exponents.
These originally appeared in \autocite{Simonenko1964}, whence they are also known as
\qu{Simonenko indices} (cf.~\autocite[20]{Maligranda1985}).
Using the notations of \autocite{BoninoCoriascoPeterssonToft2024}, we recall their definition
in the following form.


\begin{defn}
    Let \(\Phi\) be a quasi-Young function and let \(\Omega = \sets{t>0}{0<\Phi(t)<\infty}\).
    Then the \emph{Lebesgue exponents} \(p_\Phi\) and \(q_\Phi\) are given by
    \begin{align*}
        p_{\Phi}
         & \equiv
        \begin{cases}
            {\displaystyle{\sup_{t\in \Omega}
            \left (\frac{t \Phi _+ '(t)}{\Phi(t)}\right )}}, & \Omega=(0,\infty),
            \\[1ex]
            \infty,                                          & \Omega\neq (0,\infty),
        \end{cases}
        \intertext{and}
        q_{\Phi}
         & \equiv
        \begin{cases}
            {\displaystyle{\inf_{t\in \Omega}
            \left (\frac{t \Phi _+ '(t)}{\Phi(t)}\right )}}, & \Omega\neq \emptyset,
            \\
            \infty,                                          & \Omega = \emptyset.
        \end{cases}
    \end{align*}
\end{defn}


\begin{rem}
    Young functions are not necessarily differentiable, but being convex,
    they are still semi-differentiable.
    Since the definition above is the same whether one uses the left derivative or
    the right derivative, we will simply choose to use the right derivative, arbitrarily.
\end{rem}


\begin{rem}\label{rem:delta2}
    A Young function is said to fulfill the so-called \emph{\(\Delta_2\)-condition}
    (cf.~\autocite[6]{BirnbaumOrlicz1931}) if there exists a constant \(C>0\) such that
    \[
        \Phi(2t) \le C \Phi(t), \quad t\ge 0.
    \]
    It can be shown that \(\Phi\) fulfills the \(\Delta_2\)-condition if and only if
    \(p_\Phi < \infty\).
    (This is a well-known result, but for an explicit proof,
    see \autocite[Proposition 2.1]{BoninoCoriascoPeterssonToft2024}, for instance.)
\end{rem}


\begin{rem}\label{rem:pPhiquasi}
    Let \(\Phi\) be a quasi-Young function, let \(r\) be its order, and
    let \(\Psi = \Phi(t^{1/r})\), so that \(\Psi\) is a Young function.
    Then \(q_\Phi = r q_\Psi\) and \(p_\Phi = r p_\Psi\).
    Evidently, this means that \(p_\Phi < \infty\) if and only if \(p_\Psi<\infty\),
    meaning that \(\Phi\) fulfills the \(\Delta_2\)-condition if and only if
    \(\Psi\) fulfills the \(\Delta_2\)-condition.
\end{rem}

We are now equipped to recall the definitions of the various Orlicz type spaces which
we will explore in this paper.
To simplify notations,  we let
\(\rho_\Phi(f) = \nm{\Phi(|f|)}[L^1].\)


\begin{defn}
    Let \(\Phi \) be a quasi-Young function.
    The \emph{Orlicz space} \(L^{\Phi }(\rd)\) consists of all
    Lebesgue measurable functions \(\fdef\) such that
    \[
        \nm{f}[L^{\Phi}(\rd)] \equiv \nm{f}[L^\Phi] \equiv
        \inf \sets*{\lambda>0}{\rho_\Phi(f/\lambda)\leqslant 1}
    \]
    is finite.
\end{defn}


\begin{rem}\label{rem:qnorm}
    If \(\Phi\) is a quasi-Young function of order \(r\) and \(\Psi(t) = \Phi(t^{1/r})\)
    (meaning \(\Psi\) is a Young function), then
    \(\nm{\cdot}[L^\Phi] \equiv \nm{|\cdot|^r}[L^\Psi]^{1/r}\) becomes a quasi-norm of
    order \(r\) and \(L^\Phi\) a quasi-Banach space of order \(r\).
\end{rem}


\begin{rem}\label{rem:simpledense}
    We note that the analysis of \autocite[Chapter 3]{HarjulehtoHaestoe2019} can be applied in
    the case of quasi-Young functions fulfilling the \(\Delta_2\)-condition.
    In particular, Lemmas 3.1.3, 3.1.4, 3.2.4, 3.2.7, 3.2.9 and 3.2.11 carry over
    directly, as is the case for Corollary 3.2.10.
    Notably, this implies that Proposition 3.5.1 holds as well, which
    we state in this context as follows:
    if \(\Phi\) is a quasi-Young function with \(p_\Phi < \infty\), then
    the set of simple functions defined on \(\rd\) is dense in \(L^\Phi(\rd)\).
\end{rem}

\zcref{rem:simpledense} immediately gives the following known result whose proof we include
in order to be self-contained.


\begin{prop}\label{rem:Ccdense}
    If \(\Phi\) is a quasi-Young function with \(p_\Phi<\infty\), then
    \(\comp (\rd)\) is dense in \(L^\Phi(\rd)\).
\end{prop}
\begin{proof} 
    Since simple functions are dense in \(L^\Phi\) (cf.~\zcref{rem:simpledense}), it is enough
    to show that for every simple function \(f\), there exists a sequence \(f_k\in \comp\)
    such that \(\nm{f - f_k}[L^\Phi] \rightarrow 0\) whenever \(k \rightarrow \infty\).
    In fact, by linearity, it is enough to show this statement with \(f\) as
    the indicator function of a bounded measurable set.

    Let \(r\) be the order of \(\Phi\) and let \(f = \chi_A\) be the indicator function for
    a bounded measurable set \(A\subseteq \rd\).
    Let
    \[
        \fy(x) =
        \begin{cases}
            e^{-1/(1-|x|^2)} & |x|< 1,   \\
            0                & |x|\ge 1,
        \end{cases}
    \]
    \(\tilde{\fy} = \fy/\int \fy(x) \, d x\), and for each \(k\in\mb N\) let
    \(\fy_k(x) = k^d\tilde{\fy}(k x)\).
    Lastly, let \(f_k = f *\fy_{k}\).
    Then
    \[
        |f - f_k| \le g,
    \]
    where \(g = 2 \chi_B\) and
    \[
        B = \sets{x\in \rd}{|x-y|\le 1 \text{ for some }y\in A}.
    \]
    Since \(g \in L^{r}\cap L^{p_\Phi}\) and \(f_k\rightarrow f\) a.e., it follows by
    Lebesgue's dominated convergence theorem that
    \[
        \nm{f - f_k}[L^r] \rightarrow 0 \mathand \nm{f - f_k}[L^{p_\Phi}] \rightarrow 0
    \]
    whenever \(k\rightarrow \infty\).
    But since \(L^r \cap L^{p_\Phi} \subseteq L^\Phi\), there exists a constant such that
    \[
        \nm{f}[L^\Phi] \le C( \nm{f}[L^r] + \nm{f}[L^{p_\Phi}]), \quad f\in L^\Phi.
    \]
    Hence \(\nm{f - f_k}[L^\Phi] \rightarrow 0\) whenever \(k\rightarrow \infty\),
    as was to be shown.
\end{proof} 


We recall further the following two definitions.


\begin{defn}
    Let \(\Phi\) and \(\Psi\) be quasi-Young functions.
    The \emph{mixed Orlicz space} \(L^{\Phi,\Psi}(\rd[2d])\) consists of all
    Lebesgue measurable functions \(\fdef[f][\rd[2d]]\) such that
    \[
        \nm{f}[L^{\Phi,\Psi}(\rd[2d])] \equiv \nm{f}[L^{\Phi,\Psi}]\equiv \Nm{f_2^\Phi}[L^{\Psi}]
    \]
    is finite, where
    \[
        f_2^\Phi(y) = \nm{f(\cdot,y)}[L^{\Phi}].
    \]
\end{defn}


\begin{defn}\label{def:ModSpace}
    Let \(\Phi\) and \(\Psi\) be quasi-Young functions.
    The \emph{Orlicz modulation space} \(M^{\Phi,\Psi}(\rd)\) consists of all
    \(f\in\ms S'(\rd)\) such that
    \[
        \nm{f}[M^{\Phi,\Psi}(\rd)] \equiv \nm{f}[M^{\Phi,\Psi}]
        \equiv \nm{V_\fy f}[L^{\Phi,\Psi}]
    \]
    is finite, where \(\fy\in\ms S(\rd)\setminus\{0\}\) is a window function.
\end{defn}


\begin{rem} \label{rem:equivOrliczModNorm}
    With \(\tilde{V}_\fy\) as in \zcref{rem:equivModNorm}, we similarly observe that the norm
    \[
        \nm{f}[M^{\Phi,\Psi}]^*\equiv \nm{\tilde{V}_\fy f}[L^{\Phi,\Psi}],
    \]
    is equivalent to \(\nm{f}[M^{\Phi,\Psi}]\).
\end{rem}


In our investigations, we will consider \(L^p\) spaces with \(p<1\) and
\(L^\Phi\) spaces with \(\Phi\) as quasi-Young functions of order \(r<1\).
In such cases, the Fourier multiplier \(m(D)\) is not necessarily well defined, since
these spaces contain elements which are not distributions.
However, since \(m(D)\) is well defined on compactly supported functions, we can define
\(\fdef[m(D)][L^{p_1}][L^{p_2}]\), \(p_1\in(0,\infty)\), \(p_2\in(0,\infty]\) by
continuity extensions, since \(\comp\) is dense in \(L^{p_1}\).
In a similar manner, we define \(\fdef[m(D)][L^{\Phi_1}][L^{\Phi_2}]\) by continuity
extensions for quasi-Young functions \(\Phi_1\) and \(\Phi_2\).
Here, it is sufficient to assume that \(p_{\Phi_1} < \infty\), since
\(\comp\) is then dense in \(L^{\Phi_1}\) (cf.~\zcref{rem:Ccdense}).

\section{Fourier multipliers}\label{sec:fmult}
\subsection{General Orlicz space extensions}\label{subsec:general}
We begin this section with two results which we will use to
generalize results for Fourier multipliers on Lebesgue spaces to Orlicz spaces.
The first result is a special case of \autocite[Theorem 5.1]{LiuWang2012}, which
we state without proof here.


\begin{prop}\label{prop:interpol}
    Let \(q,p\in (0,\infty]\) and let \(\Phi\) be a Young function with
    \[
        q < q_\Phi \le p_\Phi <p.
    \]
    Further, let \(T\) be a linear and continuous map on \(L^{q}(\rd)+L^{p}(\rd)\) which
    restricts to linear and continuous mappings on \(L^{q}(\rd)\) and \(L^{p}(\rd)\).
    Then \(T\) is linear and continuous on \(L^{\Phi}(\rd)\) as well.
\end{prop}

The second result is a generalization of \autocite[Theorem 16]{FeichtingerNarimani2006},
whose very simple proof we present directly thereafter.


\begin{thm}\label{thm:ezproof}
    Let \(\Phi_1\) be a quasi-Young function with \(p_{\Phi_1}<\infty\) or a Young function,
    possibly with \(p_{\Phi_1}=\infty\).
    Let \(\Phi_2\) and \(\Psi\) be quasi-Young functions and suppose that
    \(\fdef[m(D)][L^{\Phi_1}(\rd)][L^{\Phi_2}(\rd)]\) is bounded.
    Then \(m(D)\) is also bounded from \(M^{\Phi_1,\Psi}(\rd)\) to \(M^{\Phi_2,\Psi}(\rd)\).
\end{thm}
\begin{proof} 
    Evidently,
    \[
        m(D_x)(\tilde{V}_\fy f)(x,\xi) = \tilde{V}_\fy( m(D) f)(x,\xi),
    \]
    and by assumption,
    \[
        \nm{m(D) g}[L^{\Phi_2}] \leq C \nm{g}[L^{\Phi_1}], \quad g\in L^{\Phi_1}(\rd),
    \]
    hence \zcref{rem:equivOrliczModNorm} gives
    \begin{align*}
        \nm{m(D)f}[M^{\Phi_2,\Psi}]
         & = \nm{\tilde{V}_\fy( m(D) f)}[L^{\Phi_2,\Psi}]  \\
         & = \nm{m(D_x)(\tilde{V}_\fy f)}[L^{\Phi_2,\Psi}] \\
         & \le C\nm{\tilde{V}_\fy f}[L^{\Phi_1,\Psi}]      \\
         & = C \nm{f}[M^{\Phi_1,\Psi}],
    \end{align*}
    which completes the proof.
\end{proof} 


\begin{rem}
    The condition \(p_{\Phi_1}<\infty\) is only included in \zcref{thm:ezproof} to
    ensure that the Fourier multiplier is well-defined on \(L^{\Phi_1}\) in the case that
    \(\Phi_1\) is a quasi-Young function of order \(r<1\).
\end{rem}

We state explicitly the following special case of \zcref{thm:ezproof},
which is a slight extension of \autocite[Theorem 16]{FeichtingerNarimani2006}.


\begin{cor}
    Let \(p_1,p_2,q \in (0,\infty]\) and suppose that
    \[
        m(D): L^{p_1}(\rd) \rightarrow L^{p_2}(\rd)
    \]
    is bounded.
    Then \(m(D)\) is bounded from \(M^{p_1,q}(\rd)\) to \(M^{p_2,q}(\rd)\) as well.
\end{cor}

We can now combine \zcref{prop:interpol,thm:ezproof} to obtain
the following extension of Mihlin's Fourier multiplier theorem.


\begin{cor}\label{cor:mihlin}
    Suppose that \(\Phi\) is a Young function with \(p_\Phi<\infty\) and \(q_\Phi > 1\),
    \(\Psi\) is a quasi-Young function, and that
    \(m\in L^\infty(\rd\setminus\{0\})\) fulfills
    \[
        \sup_{\xi\neq 0}\left(|\xi|^{|\alpha|}|\partial^\alpha m(\xi)|\right)<\infty
    \]
    for every \(\alpha\in\nd\) with \(|\alpha| \le \lfloor \frac{d}{2}\rfloor+1\).
    Then \(m(D)\) is bounded on \(M^{\Phi,\Psi}(\rd)\).
\end{cor}
\begin{proof} 
    Apply \zcref{prop:interpol,thm:ezproof} to
    Mihlin's Fourier multiplier theorem (cf.~\autocite{Mihlin1957}).
\end{proof} 

Similarly, we can extend Hörmander's improvement of
Mihlin's Fourier multiplier theorem (cf.~\autocite{Hoermander1960}).


\begin{cor}\label{cor:horm}
    Let \(\Phi\) be a Young function with \(p_\Phi<\infty\) and \(q_\Phi > 1\),
    \(\Psi\) be a quasi-Young function, and \(m\in L^\infty(\rd\setminus\{0\})\) be such that
    \[
        \sup_{R>0}\left(
        R^{-d+2|\alpha|}\int_{A_R} |\partial^\alpha m(\xi)|^2 \, d\xi\right) < \infty
    \]
    for every \(\alpha\in \nd\) with \(|\alpha| \le \lfloor \frac{d}{2}\rfloor+1\),
    where \(A_R = \sets{\xi\in\rd}{R<|\xi|<2R}\).
    Then \(m(D)\) is bounded on \(M^{\Phi,\Psi}(\rd)\).
\end{cor}


\begin{rem}
    Note that using \zcref{prop:interpol} and \zcref{thm:ezproof} to
    transfer continuity results for Fourier multipliers on Lebesgue spaces to
    (Orlicz) modulation spaces does not always lead to optimal results.
    As an example, applying \zcref{prop:interpol,thm:ezproof} to
    Theorem 7 of \autocite{FeichtingerWeisz2025} gives a result that does not cover
    Theorem 4 of the same paper.
\end{rem}


\subsection{Quasi-Banach Orlicz modulation spaces}\label{subsec:quasibanach}
We will now move on to generalizing Theorem 1 of \autocite{BenyiGroechenigOkoudjouRogers2007} to
the situation of quasi-Banach Orlicz modulation spaces.
To accomplish this, we will state and prove a series of results whose formulations mirror
those of Theorem 9, Lemma 10, and Theorem 11 of \autocite{BenyiGroechenigOkoudjouRogers2007}.
The proofs deviate to varying degrees, and most notably the proofs of
\zcref{thm:compmbound,thm:boundeddiff} below deviate significantly from
the proofs of Theorem 9 and Theorem 11, respectively.
To achieve this, we will need the following lemma about
convolution properties for Orlicz modulation spaces.


\begin{lem}\label{lem:quasiOrlConv}
    Suppose that \(\Phi\) and \(\Psi\) are quasi-Young functions and that
    \(\Phi\) is of order \(r\in(0,1]\).
    Then, for \(f\in M^{r,\infty}\) and \(g\in M^{\Phi,\Psi}\),
    \[
        \nm{f*g}[M^{\Phi,\Psi}] \leq C\nm{f}[M^{r,\infty}] \nm{g}[M^{\Phi,\Psi}]
    \]
    where \(C>0\) is a constant.
\end{lem}
\begin{proof} 
    The proof follows closely that of \autocite[Theorem 3.7]{TeofanovToft2022}.
    Let \(h = f*g\) and let \(\fy,\fy_1,\fy_2\in\ms S(\rd)\) be window functions such that
    \[
        \fy = (2\pi)^{d/2}\fy_1*\fy_2\ne 0.
    \]
    For \(j,k\in\zd\), let \(Q_r(j)\), \(Q(j)\), and \(Q(j,k)\) be as in \eqref{eq:Qj}.
    Furthermore, let
    \begin{align*}
        F(x,\xi) & = V_{\fy_1} f(x,\xi), \quad       & G(x,\xi) & = V_{\fy_2}g(x,\xi),                            \\
        H(x,\xi) & = V_\fy h(x,\xi),\qquad\text{and} & J(x,\xi) & =\left(|F(\cdot,\xi)|*|G(\cdot,\xi)|\right)(x). \\
        \intertext{Additionally, for \(j,k\in\zd\), let}
        a(j,k)   & = \nm{F}[L^{\infty}(Q(j,k))],     & b(j,k)   & = \nm{G}[L^{\infty}(Q(j,k))], \text{ and}       \\
        c(j,k)   & =\nm{J}[L^{\infty}(Q(j,k))].
    \end{align*}
    Then
    \begin{align}\label{eq:fghabm1}
        \nm{h}[M^{\Phi,\Psi}]   & = \nm{H}[L^{\Phi,\Psi}] \le \nm{H}[W^\infty(\ell^{\Phi,\Psi})]
        \le \nm{J}[W^\infty(\ell^{\Phi,\Psi})] = \nm{c}[\ell^{\Phi,\Psi}],
        \intertext{and further, there exists a constant \(C>0\) such that}
        \nm{a}[\ell^{r,\infty}] & \le C \nm{f}[M^{r,\infty}] \mathand \nm{b}[\ell^{\Phi,\Psi}]
        \le C\nm{g}[M^{\Phi,\Psi}] \label{eq:fghabm2}
    \end{align}
    (see \zcref{rem:modprop}).
    Since, for any \(x\in \rd\) and \(n \in \mb Z^d\),
    \[
        |F(x,n)|\le \sum_{j\in\mb Z^d} a(j,n)\chi_{Q(j)}(x) \quad \text{and} \quad
        |G(x,n)|\le \sum_{k\in\mb Z^d} b(k,n)\chi_{Q(k)}(x),
    \]
    it follows that

    \begin{align*}
        J(x,n)
         & \le \sum_{j,k\in\mb Z^d} a(j,n) b(k,n) (\chi_{Q(j)}*\chi_{Q(k)})(x) \\
         & \le \sum_{j,k\in\mb Z^d} a(j,n) b(k,n) \chi_{Q_2(j+k)}(x).
    \end{align*}
    For any \(l\in \mb Z^d\), let
    \[
        \Omega_l = \sets{(j,k)\in\mb Z^{2d}}{j+k\in[l-2v_1,l]^d},
    \]
    where \(v_1=(1,1,\dots,1)\).
    Then \(l\in Q_2(j+k)\) if and only if \((j,k)\in \Omega_l\), hence
    \begin{align*}
        \sum_{j,k\in\mb Z^d} a(j,n) b(k,n) \chi_{A_{j,k}}(l)
         & = \sum_{(j,k)\in\Omega_l}a(j,n) b(k,n)                              \\
         & = \sum_{t\le 2v_1} \big( a(\cdot,n) * b(\cdot,n) \big)(t+l - 2v_1).
    \end{align*}
    Since this sum is finite, we obtain for some constants \(C_j>0\), \(j=0,1,2\),
    \begin{align*}
        \nm{c(\cdot,n)}[\ell^\Phi]
         & \le \nm*{\sum_{t\le 2v_1} \big( a(\cdot,n) * b(\cdot,n) \big)(t+\cdot - 2v_1)}[\ell^\Phi]     \\
         & \le C_0 \sum_{t\le 2v_1} \nm*{\big( a(\cdot,n) * b(\cdot,n) \big)(t+\cdot - 2v_1)}[\ell^\Phi] \\
         & \le C_1 \nm*{a(\cdot,n)*b(\cdot,n)}[\ell^\Phi]                                                \\
         & \le C_2 \nm{a(\cdot,n)}[\ell^r] \nm{b(\cdot,n)}[\ell^\Phi],
    \end{align*}
    where we refer to \autocite[Lemma 4.1]{ToftUesterNabizadehOeztop2022} for the fourth inequality.
    Therefore
    \begin{align*}
        \nm{c}[\ell^{\Phi,\Psi}]
         & \le C_2 \Nm{\nm{a(\cdot,n)}[\ell^r] \nm{b(\cdot,n)}[\ell^\Phi]}[\ell^\Psi] \\
         & \le C_2 \nm{a}[\ell^{r,\infty}] \nm{b}[\ell^{\Phi,\Psi}].
    \end{align*}
    By combining \eqref{eq:fghabm1} and \eqref{eq:fghabm2} we therefore obtain
    \begin{equation*}
        \nm{h}[M^{\Phi,\Psi}] \le C_3 \nm{f}[M^{r,\infty}] \nm{g}[M^{\Phi,\Psi}]
    \end{equation*}
    for some constant \(C_3>0\), which is the desired result.
\end{proof} 


From this, we immediately obtain the following.


\begin{thm}\label{prop:MPhiPsiBound}
    Suppose that \(\Phi\) and \(\Psi\) are quasi-Young functions such that
    \(\Phi\) is of order \(r\in(0,1]\), and \(m\in W^{\infty,r}(\rd)\).
    Then \(m(D)\) is bounded on \(M^{\Phi,\Psi}(\rd)\).
\end{thm}
\begin{proof} 
    Since
    \[
        m(D)f(x) = (\ms F^{-1}m*f)(x),
    \]
    and \(m\in W^{\infty,r}\) if and only if \(\ms F^{-1}m \in M^{r,\infty}\),
    the result immediately follows from \zcref{lem:quasiOrlConv}.
\end{proof} 


We state explicitly the following special case.


\begin{cor}
    Suppose that \(\Phi,\Psi\) are Young functions and \(m\in W^{\infty,1}(\rd)\).
    Then \(m(D)\) is bounded on \(M^{\Phi,\Psi}(\rd)\).
\end{cor}

Before we state the next theorem in the sequence of results mirroring those of
\autocite{BenyiGroechenigOkoudjouRogers2007}, we will need the following lemma, which is a
direct consequence of \autocite[Theorem 3.3]{Toft2022}.


\begin{lem}\label{lem:compDist}
    Let \(f\in \ms E'(\rd)\) and let \(K=\supp f\). Then, for \(0<p_1,p_2,q\le\infty\),
    \[
        C_K^{-1}\nm{f}[W^{p_1,q}] \le \nm{f}[W^{p_2,q}] \le C_K \nm{f}[W^{p_1,q}],
    \]
    where \(C_K\) is a constant depending on \(K\).
\end{lem}
\begin{proof} 
    Without loss of generality, suppose that \(p_1\le p_2\),
    let \(r = \frac{1}{p_1} - \frac{1}{p_2} > 0\), and \(s=\max\{1,q\}\).
    Since \(f\) has compact support, there exists \(\fy \in \comp\) with \(\fy = 1\) on
    \(K\) and \(\nm{\fy}[W^{r,s}] < \infty\).
    By \autocite[Theorem 3.3]{Toft2022}, there exists a constant \(C_1 > 0\) such that
    \[
        \nm{f}[W^{p_1,q}] = \nm{\fy \cdot f}[W^{p_1,q}]
        \le C_1 \nm{\fy}[W^{r,s}] \nm{f}[W^{p_2,q}] = C_2 \nm{f}[W^{p_2,q}].
    \]
    This completes the proof, since
    \(\nm{f}[W^{p_2,q}] \le C_3 \nm{f}[W^{p_1,q}]\) trivially holds when
    \(p_1 \le p_2\) for some constant \(C_3>0\).
\end{proof} 


With this result in mind, we prove the following generalization of
\autocite[Theorem 9]{BenyiGroechenigOkoudjouRogers2007}.
As mentioned before, the conditions in the statement of the theorem need only be
slightly altered, but the strategy of the proof differs.


\begin{thm}\label{thm:compmbound}
    Suppose that \(\Phi\) and \(\Psi\) are quasi-Young functions,
    where \(\Phi\) is of order \(r\in(0,1]\), and that
    \(\mu\in C^{N}(\rd\setminus\{0\})\), where \(N>\frac{d}{r}\) is an integer and
    \(\mu\) is homogeneous of order \(\alpha> \frac{d(1-r)}{r}\).
    Let \(\chi\in\comp(\rd;[0,1])\) fulfill
    \begin{align*}
        \chi(\xi) =
        \begin{cases}
            1, & |\xi|\le 1, \\
            0, & |\xi|\ge 2,
        \end{cases}
    \end{align*}
    and let \(m = e^{i\mu}\chi\).
    Then the following holds:
    \begin{enumerate}
        \item \label{thm:compbounditem1} \(m\in \bigcap_{p>0} W^{p,r}(\rd)\);
              \vrum
        \item \label{thm:compbounditem2} \(m(D)\) is bounded on \(M^{\Phi,\Psi}(\rd)\).
    \end{enumerate}
\end{thm}
\begin{proof} 
    We begin the proof in the same vein as in that of \autocite[Theorem 9]{BenyiGroechenigOkoudjouRogers2007}.
    By \zcref{prop:MPhiPsiBound},
    \ref{thm:compbounditem2} follows from \ref{thm:compbounditem1}.
    Since \(m\) has compact support, it is enough to show that
    \(m\in W^{\infty,r}\) by \zcref{lem:compDist}.

    By Taylor expansion,
    \begin{equation}\label{eq:mtaylor}
        m(\xi) = \sum_{k=0}^\infty \frac{i^k}{k!} \phi_k(\xi),
    \end{equation}
    with \(\phi_k(\xi) = \mu^k(\xi)\chi(\xi)\).
    Letting \(\psi(\xi) = \chi(\xi/2) - \chi(\xi)\), we have
    \[
        \chi(\xi) =\sum_{j=1}^\infty \psi(2^j\xi),
    \]
    so that, by the homogeneity of \(\mu\),
    \[
        \phi_k(\xi) = \sum_{j=1}^\infty \mu^k(\xi)\psi(2^j \xi)
        = \sum_{j=1}^\infty 2^{-j k \alpha} \psi_k(2^j \xi),
    \]
    where, in turn, \(\psi_k(\xi) = \mu^k(\xi)\psi(\xi)\).

    Now, take any window function \(\fy\in \ms S\setminus\{0\}\)
    (cf.~\autocite[Theorem 3.1]{GalperinSamarah2004}).
    Then
    \begin{align}\label{eq:spliteq}
        \begin{split}
            \nm{\phi_k}[W^{\infty,r}]^r
             & \le \sum_{j=1}^\infty 2^{-jk\alpha r} \nm{\psi_k(2^j\cdot)}[W^{\infty,r}]^r                    \\
             & = \sum_{j=1}^\infty 2^{-jk\alpha r} \nm{V_\fy \left(\psi_k(2^j\cdot)\right)}[L_*^{\infty,r}]^r \\
             & = \sum_{j=1}^\infty 2^{-jk\alpha r+jd(1-r)}\nm{V_{\fy_j} \psi_k}[L_*^{\infty,r}]^r,
        \end{split}
    \end{align}
    where \(\fy_j = \fy(2^{-j}\cdot)\).
    Let \(\tilde{\fy}\in\comp(\rd;[0,1])\) fulfill
    \[
        \tilde{\fy}(\xi) =
        \begin{cases}
            0, & |\xi|<\frac12 \text{ or } |\xi|\ge5, \\
            1, & 1\le |\xi| \le 4.
        \end{cases}
    \]
    Then, since \(\tilde{\fy} = 1\) on \(\supp \psi_k\),
    \begin{align*}
        V_{\fy_j} \psi_k(\xi,x) & = (2\pi)^{-d/2} \ms F\left[\psi_k \fy_j(\cdot-\xi) \right](x)                                                     \\
                                & = (2\pi)^{-d/2}\ms F[\psi_k \tilde{\fy}^2 \fy_j(\cdot-\xi)](x)                                                    \\
                                & = (2\pi)^{-d/2}\Big(\ms F\left[\psi_k \tilde{\fy} \right]* \ms F\left[\tilde{\fy} \fy_j(\cdot-\xi)\right]\Big)(x) \\
                                & = (2\pi)^{-d/2}\Big(V_{\tilde{\fy}}\psi_k(0,\cdot) * V_{\fy_j} \tilde{\fy}(\xi,\cdot) \Big)(x).
    \end{align*}
    Now, for any \(n=1,\dots,d\), and any integer \(0\le N_0\le N\),
    \begin{align*}
        (2\pi)^{d/2} \left|x_n^{N_0} V_{\fy_j} \tilde{\fy}(\xi,x) \right|
         & = \left|\int \tilde{\fy}(\eta) \fy_j(\eta-\xi)\partial_n^{N_0} e^{-i \scal{\eta}{x}} \, d \eta \right|                          \\
         & \le C_{N} \sum_{N_1+N_2 = N_0} \int \left|(\partial_n^{N_1}\tilde{\fy})(\eta)(\partial_n^{N_2}\fy_j)(\eta-\xi) \right| \, d\eta \\
         & \le \tilde{C}_{N},
    \end{align*}
    where \(\tilde{C}_{N}\) is a constant depending on \(N\) and on
    \[
        \sup_{N_1 + N_2 \le N} \nm{\partial_n^{N_1} \tilde{\fy}}[L^1]
        \nm{\partial_n^{N_2} \fy}[L^\infty] < \infty.
    \]
    Hence
    \begin{align*}
        \nm{V_{\fy_j}\psi_k}[L_*^{\infty,r}]^r
         & \le C_N \nm{\left(V_{\tilde{\fy}}\psi_k(0,\cdot) *G\right)(x)}[L^r]^r,
    \end{align*}
    for some new constant \(C_N\) only depending on \(N\), where
    \begin{equation}\label{eq:G(y)}
        G(y) = \left( \max_{1\le n \le d}\left\{ |y_n|^N,1\right\} \right)^{-1}.
    \end{equation}
    We let \(Q(j)\) be as in the proof of \zcref{lem:quasiOrlConv},
    but this time we let
    \begin{align*}
        F(x) & = V_{\tilde{\fy}}\psi_k(0,x), & J(x) & = \left(|F|*|G|\right)(x),             \\
        \intertext{and}
        a(l) & = \nm{F}[L^{\infty}(Q(l))],   & b(l) & = \nm{G}[L^{\infty}(Q(l))],\text{ and} \\
        c(l) & = \Nm{J}[L^{\infty}(Q(l))].
    \end{align*}
    Then
    \[
        |F(x)| \le \sum_{l\in\zd} a(l) \chi_{Q(l)}(x) \quad \text{and} \quad
        |G(x)| \le \sum_{l\in\zd} b(l) \chi_{Q(l)}(x),
    \]
    and by proceeding analogously to the proof of \zcref{lem:quasiOrlConv},
    we obtain (with \(\Phi(t) = t^r\))
    \[
        \nm{V_{\fy_j}\psi_k}[L_*^{\infty,r}] \le C_N \nm{J}[L^r]^r
        \le C_N\nm{c}[\ell^r]^r \le C \nm{a}[\ell^r]^r \nm{b}[\ell^r]^r,
    \]
    for some constant \(C>0\).
    Since \(N>\frac{d}{r}\),
    \[
        \nm{b}[\ell^r]^r = \sum_{l\in\zd} G(l)^r < \infty,
    \]
    and for \(\nm{a}[\ell^r]^r\) we have
    \begin{align*}
        \nm{a}[\ell^r]^r
         & = \sum_{l\in \zd} \nm{F}[L^\infty(Q(l))]^r
        \le \sum_{l\in\zd} \nm{V_{\tilde{\fy}} \psi_k}[L^\infty(Q(0,l))]^r                                \\
         & \le \sup_{\lambda\in\rd} \sum_{l\in \zd} \nm{V_{\tilde{\fy}} \psi_k}[L^\infty(Q(\lambda,l))]^r
        = \nm{V_{\tilde{\fy}} \psi_k}[W^\infty(\ell_*^{\infty,r})]^r                                      \\
         & \le C\nm{V_{\tilde{\fy}}\psi_k}[L_*^{\infty,r}]^r,
    \end{align*}
    where we refer to \autocite[Theorem 3.3]{GalperinSamarah2004} in the last inequality.
    Hence, by \eqref{eq:spliteq},
    \[
        \nm{\phi_k}[W^{\infty,r}]^r
        \le  C_{r,N} \nm{V_{\tilde{\fy}}\psi_k}[L_*^{\infty,r}]^r
        \sum_{j=1}^\infty 2^{-j k \alpha r + j d (1 - r)}
        \le C_{r,N} 2^{-\beta_k}\nm{V_{\tilde{\fy}} \psi_k}[L_*^{\infty,r}]^r,
    \]
    where \(\beta_k = k\alpha r - d(1-r) >0\) by assumption.

    It remains to find an upper bound for \(\nm{V_{\tilde{\fy}} \psi_k}[L_*^{\infty,r}]\) so that
    the estimate of \(\nm{\phi_k}[W^{\infty,r}]^r\) from \eqref{eq:spliteq} together with
    \eqref{eq:mtaylor} gives \(\nm{m}[W^{\infty,r}]<\infty\), which is the desired result.

    We have, for any \(n=1,\dots,d\),
    (noting that \(\supp \psi \subseteq \sets{\eta\in\rd}{1\le |\eta|\le 4}\))
    \begin{align*}
         & \left| x_n^N V_{\tilde{\fy}} \psi_k(\xi,x) \right|                                                       \\
         & = (2\pi)^{-d/2}
        \left|\int \partial_n^N \big( \psi_k(\eta) \tilde{\fy}(\eta-\xi) \big) e^{-i \scal{\eta}{x}} d \eta \right| \\
         & \le (2\pi)^{-d/2} \sum_{|\gamma|+N_1+N_2 = N} \frac{N!}{\gamma!N_1! N_2!}
        \int \left( \prod_{j=1}^k |\partial_n^{\gamma_j} \mu(\eta)| \right)
        |\partial_n^{N_1}\psi(\eta)| |\partial_n^{N_2}\tilde{\fy}(\eta-\xi)| d \eta                                 \\
         & \le C(\xi) \sum_{|\gamma|+N_1+N_2 = N} \frac{N!}{\gamma!N_1! N_2!}
        \prod_{j=1}^k \left( C_1 4^{\alpha - |\gamma_j|}\right)                                                     \\
         & \le C(\xi) C_2^k 4^{k \alpha},
    \end{align*}
    where \(\gamma=(\gamma_1,\dots,\gamma_k)\in \mb N^k\),
    \[
        C_1 =
        \max_{\gamma_j\le N} \left( \sup_{|\eta|=1}|\partial_n^{\gamma_j}\mu(\eta)| \right),
    \]
    and where we use the homogeneity of \(\mu\) for the second inequality.
    Hence
    \[
        |V_{\tilde{\fy}} \psi_k (\xi,x)| \le C(\xi) C_2^k 4^{k \alpha} G(x).
    \]
    Finally, since \(N > \frac{d}{r}\), this gives us
    \[
        \nm{V_{\tilde{\fy}}\psi_k}[L_*^{\infty,r}]^r
        \le C_3 C_2^{kr} 4^{\alpha k r} \nm{G}[L^r]^r \le C_4 C_2^{k r} 4^{\alpha k r}.
    \]
    Using this with \eqref{eq:mtaylor} and \eqref{eq:spliteq}, we now get (since \(r\le 1\))
    \begin{align*}
        \nm{m}[W^{\infty,r}]^r
         & \le \sum_{k=0}^\infty \frac{1}{(k!)^r} \nm{\phi_k}[W^{\infty,r}]^r \\
         & \le \sum_{k=0}^\infty\frac{2^{-\beta_k}}{(k!)^r}
        \nm{V_{\tilde{\fy}}\psi_k}[L_*^{\infty,r}]^r                          \\
         & \le C_4 2^{d(1-r)} \sum_{k=0}^\infty
        \frac{\left(C_2^r 2^{\alpha r} \right)^k}{(k!)^r} < \infty,
    \end{align*}
    which, at last, gives the desired result.
\end{proof} 


To obtain the final result needed in the previously mentioned sequence of results from
\autocite{BenyiGroechenigOkoudjouRogers2007}, we will make use of the following lemma, which is almost
completely identical to \autocite[Lemma 10]{BenyiGroechenigOkoudjouRogers2007}.
We include its short proof for the sake of completeness.


\begin{lem}\label{lem:mtilde}
    Suppose that \(\fdef[f][\rd][\mb R]\) and \(\fdef[g][\rd][\rd]\) are
    measurable functions and let
    \[
        \tilde{m}_\xi(\eta) = m(\eta) e^{i(f(\xi)+\scal{\eta}{g(\xi)})}.
    \]
    Then
    \[
        \nm{m}[W^{\infty,r}] = \sup_{\xi\in\rd} \nm{V_\fy \tilde{m}_\xi(\xi,\cdot)}[L^r],
    \]
    for some window function \(\fy\in\ms S\).
\end{lem}
\begin{proof} 
    We have
    \begin{align*}
        \sup_{\xi\in\rd} \nm{V_\fy \tilde{m}_\xi(\xi,\cdot)}[L^r]^r
         & = \sup_{\xi\in\rd} \int \left| (2\pi)^{-d/2}
        \int m(\eta)e^{i(f(\xi)+\scal{\eta}{g(\xi)})} \fy(\eta - \xi)
        e^{-i\scal{\eta}{x}} d\eta \right|^r d x                         \\
         & = \sup_{\xi\in\rd} \int \left| (2\pi)^{-d/2}
        \int m(\eta) \fy(\eta - \xi)
        e^{-i\scal{\eta}{x - g(\xi)}}  d\eta \right|^r d x               \\
         & = \sup_{\xi \in \rd} \nm{V_\fy m(\xi, \cdot - g(\xi))}[L^r]^r \\
         & = \sup_{\xi \in \rd} \nm{V_\fy m(\xi, \cdot)}[L^r]^r
        = \nm{m}[W^{\infty,r}]^r.
    \end{align*}
\end{proof} 


We can now obtain the following generalization of \autocite[Theorem 11]{BenyiGroechenigOkoudjouRogers2007}.


\begin{thm}\label{thm:boundeddiff}
    Let \(N >\frac{d}{r}\) be an integer. Suppose that \(\mu\in \cont{N}(\rd)\) is real-valued and fulfills
    \[
        \nm{\partial^\beta\mu}[L^\infty] \le C, \quad 2\le |\beta| \le N,
    \]
    for some constant \(C>0\), and let \(m = e^{i \mu}\). Then the following holds:
    \begin{enumerate}
        \item \(m\in W^{\infty,r}(\rd)\);
              \vrum
        \item \(m(D)\) is bounded on \(M^{\Phi,\Psi}(\rd)\) for any quasi-Young functions
              \(\Phi, \Psi\) with \(\Phi\) of order \(r\in(0,1]\).
    \end{enumerate}
\end{thm}
\begin{proof} 
    As in the proof of \zcref{thm:compmbound}, we show the first statement only, since
    it implies the second one.

    Let
    \[
        r_\xi(\eta) = \mu(\eta) - \mu(\xi) - \scal{\nabla \mu(\xi)}{\eta-\xi},
    \]
    so that by Taylor's theorem, for every \(j=1,\dots,d\),
    \begin{equation}\label{eq:r_xi(eta)2}
        |\partial_j r_\xi(\eta)| \le C |\eta - \xi|, \quad \eta\in \rd,
    \end{equation}
    and
    \begin{equation}\label{eq:r_xi(eta)3}
        \nm{\partial_j^k r_\xi}[L^\infty] \le C, \quad k=2,\dots,N.
    \end{equation}
    Let
    \[
        F(\xi,x) = (2\pi)^{-d/2} \int e^{i r_\xi(\eta)} \overline{\fy(\eta-\xi)} e^{-i\scal{\eta}{x}} d\eta
        = V_\fy \tilde m_\xi(\xi,x),
    \]
    where \(\tilde m_\xi(\eta) = e^{i f(\xi) + \scal{\eta}{g(\xi)}}\) for
    \(f(\xi)= \scal{\nabla \mu(\xi)}{\xi}-\mu(\xi)\) and \(g(\xi)=\nabla \mu(\xi)\).
    Then for \(j=1,\dots,d\) and any integer \(0\le N_0 \le N\),
    \begin{align*}
        \left|x_j^{N_0} F(\xi,x) \right|
         & = \left|\int \partial_j^{N_0} \left( e^{i r_\xi(\eta)} \overline{\fy(\eta-\xi)} \right)
        e^{-i\scal{\eta}{x}} d\eta\right|                                                          \\
         & \le C_{N_0}\sum_{k=0}^{N_0}
        \int \left| \partial_j^{k} e^{i r_\xi(\eta)} \right|
        \left| \partial_j^{N_0-k} \fy(\eta-\xi) \right| d\eta.
    \end{align*}
    Let \(T_{j,k}(f) = (\partial_j f, \partial_j^2 f, \dots, \partial_j^{k}f)\).
    Then
    \[
        |\partial_j^{k} e^{i r_\xi(\eta)}| = |P_{k}\circ T_{j,k}(r_\xi(\eta))|,
    \]
    where \(P_k\) is some polynomial of order \(k\).
    Hence, by 
    \eqref{eq:r_xi(eta)2}, \eqref{eq:r_xi(eta)3}, and the fact that \(\fy\in\ms S\),
    \begin{align*}
        \left|x_j^{N_0} F(\xi,x) \right|
         & \le C'_{N_0} \sup_{ 0\le n \le N_0}
        \left(\int |\eta-\xi|^{N_0} \partial_j^n\fy(\eta-\xi) \, d\eta \right)
        = C < \infty.
    \end{align*}
    for some (new) constants \(C'_{N_0}, C>0\).

    As in the proof of \zcref{thm:compmbound}, let \(G\) be given by \eqref{eq:G(y)}.
    The calculations above give
    \[
        \left| F(\xi,x)\right| \le C G(x),
    \]
    for every \(\xi,x\in\rd\).
    Since \(N> \frac{d}{r}\), we have \(\nm{G}[L^r]<\infty\).
    Therefore, by \zcref{lem:mtilde},
    \begin{align*}
        \nm{m}[W^{\infty,r}] = \nm{F}[L_*^{\infty,r}] \le C\nm{G}[L^r] < \infty,
    \end{align*}
    and the result follows.
\end{proof} 


This leaves us with all the tools necessary to prove the following generalization of
\autocite[Theorem 1]{BenyiGroechenigOkoudjouRogers2007} to the case of quasi-Banach Orlicz modulation spaces.


\begin{thm}\label{thm:thm1benyi}
    Let \(r\in\left(\frac{d}{d+2},1\right]\), \(\alpha\in\left(\frac{d(1-r)}{r}, 2\right]\),
    and \(m(\xi) = e^{i \mu(\xi)}\), where \(\mu\in C^N(\rd\setminus\{0\})\) is real-valued and
    homogeneous of order \(\alpha\), and where \(N>\frac{d}{r}\) is an integer.
    Then the following holds:
    \begin{enumerate}
        \item \(m\in W^{\infty,r}(\rd)\);
              \vrum
        \item \(m(D)\) is bounded on \(M^{\Phi,\Psi}(\rd)\) for any quasi-Young functions
              \(\Phi, \Psi\), with \(\Phi\) of order \(r\).
    \end{enumerate}
\end{thm}
\begin{proof} 
    The condition that \(r\in\left(\frac{d}{d+2},1\right]\) ensures that \(\frac{d(1-r)}{r}<2\),
    so that there exists \(\alpha\) fulfilling both
    \(\alpha > \frac{d(1-r)}{r}\) and \(\alpha \le 2\).

    Let \(\chi\) be as in \zcref{thm:compmbound},
    \[
        m_1(\xi) = e^{i \mu(\xi)} \chi(\xi),\mathandd m_2(\xi) = e^{i\mu(\xi)}(1-\chi(\xi)).
    \]
    Then \(m=m_1+m_2\).
    Since \(m_1\in W^{\infty,r}\) by \zcref{thm:compmbound} (\(\alpha > \frac{d(1-r)}{r}\)), and
    \(m_2\in W^{\infty,r}\) by \zcref{thm:boundeddiff} (\(\alpha \le 2\)), it follows that
    \(m\in W^{\infty,r}\).
    Hence, the proof is complete by \zcref{prop:MPhiPsiBound}.
\end{proof} 

We state explicitly the following special case.


\begin{cor}\label{cor:r1}
    Let \(\Phi\), \(\Psi\) be Young functions, \(\alpha\in(0,2]\), \(
    a \in \mb R\), and
    \(m(\xi) = e^{i a |\xi|^\alpha}\).
    Then \(m(D)\) is bounded on \(M^{\Phi,\Psi}(\rd)\).
\end{cor}

If \(\alpha =2\), we can remove the restrictions on the order \(r\) of
the quasi-Young function \(\Phi\) in \zcref{thm:thm1benyi}.


\begin{thm}\label{thm:alpha2}
    Let \(A\) by a real \(d\times d\) matrix and let \(m(\xi) = e^{i \scal{A \xi}{\xi}}\).
    Then \(m(D)\) is a bounded operator on \(M^{\Phi,\Psi}\), where
    \(\Phi\) and \(\Psi\) are quasi-Young functions.
\end{thm}
\begin{proof} 
    Suppose that \(\fy\in \ms S\setminus\{0\}\) and \(f\in M^{\Phi,\Psi}\).
    For some \(d \times d\) matrix \(B_A\),
    \[
        \left|V_\fy \left( e^{i \scal{A D}{D}} f\right)(x,\xi) \right|
        = |V_{\fy_A} f(x+B_A\xi, \xi)| \equiv |T_A f|,
    \]
    where \(\fy_A = e^{-i\scal{A D}{D}}\fy \in \ms S\setminus\{0\}\)
    (cf.~\autocite[Prop. 1.5]{Toft2004a}).
    By the translation invariance of \(L^\Phi\) and by \zcref{rem:modprop},
    \[
        \nm{T_{A} f }[L^{\Phi,\Psi}]
        =       \nm{V_{\fy_A} f}[L^{\Phi,\Psi}]
        \le   C \nm{V_\fy f}[L^{\Phi,\Psi}]
        =     C \nm{f}[M^{\Phi,\Psi}]
    \]
    for some constant \(C>0\), which completes the proof.
\end{proof} 


We conclude this paper with an application of \zcref{thm:thm1benyi,thm:alpha2,cor:r1}.
It is the so-called space-fractional Schrödinger equation \autocite{Laskin2002},
which has uses in quantum mechanics and optics \autocite{Longhi2015}.


\begin{example}[name=Fractional Schrödinger equation]
    Let \(r\in\left(\frac{d}{d+2},1\right]\), \(\alpha\in\left(\frac{d(1-r)}{r}, 2\right]\),
    and \(a\in \mb R\), and consider the initial value problem
    \[
        \left\{
        \begin{aligned}
            \partial_t \psi(t,x) & = i a (-\Delta_x)^{\alpha/2} \psi(t,x), \\
            \psi(0,x)            & = \psi_0(x)\in M^{\Phi,\Psi}(\rd),
        \end{aligned}\right.
    \]
    for some quasi-Young functions  \(\Phi\), \(\Psi\) with \(\Phi\) of order \(r\).
    The solution operator is a Fourier multiplier with symbol \(e^{i a|\xi|^\alpha}\), hence
    by \zcref{thm:thm1benyi}, the solution \(\psi(t,x)\) fulfills
    \[
        \nm{\psi(t,\cdot)}[M^{\Phi,\Psi}] \le C \nm{\psi_0}[M^{\Phi,\Psi}]
    \]
    for some constant \(C>0\) depending on \(a\in \mb C\) and \(t\).

    The choices of \(r\) affect the possible choices of \(\alpha\), and vice versa.
    If \(\alpha = 2\), then we are in the standard Schrödinger case and the statement will be
    true for any \(r>0\) (as in \zcref{thm:alpha2}).
    If \(r=1\), then we are in a Banach setting, and the statement is true for any
    \(\alpha \in (0,2]\) (as in \zcref{cor:r1}).
    The special case of \(d/(d+2)<r< 1\) still allows \(\alpha < 2\), provided
    \(\alpha > d(1-r)/r\) (by \zcref{thm:thm1benyi}).
\end{example}

\section*{Acknowledgements}
The author has been supported by Vetenskapsrådet (Swedish Science Council)
within the project 2019-04890.
\printbibliography

@Article{Aoki1942,
  author       = {Aoki, Tosio},
  title        = {Locally bounded linear topological spaces},
  issn         = {0386-2194},
  number       = {10},
  pages        = {588--594},
  volume       = {18},
  date         = {1942-01},
  doi          = {10.3792/pia/1195573733},

  journaltitle = {Proc. Jpn. Acad. A: Math. Sci.},
  publisher    = {Project Euclid},
}

@Article{BenyiGroechenigOkoudjouRogers2007,
  author       = {B{\'e}nyi, {\'A}rp{\'a}d and Gr{\"o}chenig, Karlheinz and Okoudjou, Kasso A. and Rogers, Luke G.},
  title        = {Unimodular {F}ourier multipliers for Modulation Spaces},
  issn         = {0022-1236},
  number       = {2},
  pages        = {366--384},
  volume       = {246},
  date         = {2007-05},
  doi          = {10.1016/j.jfa.2006.12.019},

  journaltitle = {J. Funct. Anal.},
  publisher    = {Elsevier BV},
}

@Article{BirnbaumOrlicz1931,
  author       = {Birnbaum, Z. and Orlicz, W.},
  title        = {{\"U}ber die {V}erallgemeinerung des {B}egriffes der zueinander konjugierten {P}otenzen},
  issn         = {1730-6337},
  number       = {1},
  pages        = {1--67},
  volume       = {3},
  date         = {1931},
  doi          = {10.4064/sm-3-1-1-67},

  journaltitle = {Stud. Math.},
  language     = {German},
  publisher    = {Institute of Mathematics, Polish Academy of Sciences},
}

@Article{BoninoCoriascoPeterssonToft2024,
  author       = {Bonino, Matteo and Coriasco, Sandro and Petersson, Albin and Toft, Joachim},
  title        = {Fourier Type Operators on {O}rlicz Spaces and the Role of {O}rlicz {L}ebesgue Exponents},
  issn         = {1660-5454},
  note         = {art. no. 219},
  volume       = {21},
  date         = {2024-11-21},
  doi          = {10.1007/s00009-024-02735-9},
  journaltitle = {Mediterr. J. Math.},
  publisher    = {Springer Science and Business Media LLC},
}

@Article{FeichtingerNarimani2006,
  author       = {Hans G. Feichtinger and Ghassem Narimani},
  title        = {Fourier Multipliers of Classical Modulation Spaces},
  issn         = {1063-5203},
  number       = {3},
  pages        = {349--359},
  volume       = {21},
  date         = {2006-11},
  doi          = {10.1016/j.acha.2006.04.010},

  journaltitle = {Appl. Comput. Harmon. A.},
  publisher    = {Elsevier BV},
}

@Article{FeichtingerWeisz2025,
  author       = {Feichtinger, Hans G. and Weisz, Ferenc},
  title        = {Piecewise linear and step Fourier multipliers for modulation spaces},
  issn         = {0022-1236},
  note         = {art. no. 110795},
  number       = {5},
  volume       = {288},
  date         = {2025-03},
  doi          = {10.1016/j.jfa.2024.110795},

  journaltitle = {J. Funct. Anal.},
  publisher    = {Elsevier BV},
}

@Article{GalperinSamarah2004,
  author       = {Galperin, Yevgeniy V. and Samarah, Salti},
  title        = {Time-frequency analysis on modulation spaces $M^{p,q}_m$, $0<p,q\leq \infty$},
  issn         = {1063-5203},
  number       = {1},
  pages        = {1--18},
  volume       = {16},
  date         = {2004-01},
  doi          = {10.1016/j.acha.2003.09.001},

  journaltitle = {Appl. Comput. Harmon. A.},
  publisher    = {Elsevier BV},
}

@Article{GumberRanaToftUester2024,
  author       = {Gumber, Anupam and Rana, Nimit and Toft, Joachim and {\"U}ster, R{\"u}ya},
  title        = {Pseudo-differential Calculi and Entropy Estimates with {O}rlicz Modulation Spaces},
  issn         = {0022-1236},
  note         = {art. no. 110225},
  number       = {3},
  volume       = {286},
  date         = {2024-02},
  doi          = {10.1016/j.jfa.2023.110225},

  journaltitle = {J. Funct. Anal.},
  publisher    = {Elsevier BV},
}

@Book{HarjulehtoHaestoe2019,
  author       = {Harjulehto, Petteri and H{\"a}st{\"o}, Peter},
  publisher    = {Springer International Publishing},
  title        = {{O}rlicz Spaces and Generalized {O}rlicz Spaces},
  isbn         = {9783030151003},
  date         = {2019},
  doi          = {10.1007/978-3-030-15100-3},

  fjournal     = {Lecture Notes in Mathematics},
  issn         = {1617-9692},
  journaltitle = {Lect. Notes Math.},
}

@Article{Hoermander1960,
  author       = {H{\"o}rmander, Lars},
  title        = {Estimates for Translation Invariant Operators in $L^p$ Spaces},
  issn         = {0001-5962},
  number       = {1–2},
  pages        = {93--140},
  volume       = {104},
  date         = {1960},
  doi          = {10.1007/bf02547187},

  journaltitle = {Acta Math.},
  publisher    = {International Press of Boston},
}

@Article{Laskin2002,
  author    = {Laskin, Nick},
  journal   = {Phys. Rev. E},
  title     = {Fractional {Schrödinger} equation},
  year      = {2002},
  issn      = {1095-3787},
  month     = Nov,
  number    = {5},
  pages     = {056108},
  volume    = {66},
  doi       = {10.1103/physreve.66.056108},
  publisher = {American Physical Society (APS)},
}

@Article{LiebSolovej1991,
  author       = {Lieb, Elliott H. and Solovej, Jan Philip},
  title        = {Quantum coherent operators: A generalization of coherent states},
  issn         = {1573-0530},
  number       = {2},
  pages        = {145--154},
  volume       = {22},
  date         = {1991-06},
  doi          = {10.1007/bf00405179},
  journaltitle = {Lett. Math. Phys.},
  publisher    = {Springer Science and Business Media LLC},
}

@Article{LiuWang2012,
  author       = {Liu, PeiDe and Wang, MaoFa},
  title        = {Weak {O}rlicz Spaces: some Basic Properties and their Applications to Harmonic Analysis},
  issn         = {1869-1862},
  number       = {4},
  pages        = {789--802},
  volume       = {56},
  date         = {2012-07},
  doi          = {10.1007/s11425-012-4452-5},

  journaltitle = {Sci. China Math.},
  publisher    = {Springer Science and Business Media LLC},
}

@Article{Longhi2015,
  author    = {Longhi, Stefano},
  journal   = {Opt. Lett.},
  title     = {Fractional {Schrödinger} equation in optics},
  year      = {2015},
  issn      = {1539-4794},
  month     = Mar,
  number    = {6},
  pages     = {1117},
  volume    = {40},
  doi       = {10.1364/ol.40.001117},
  publisher = {Optica Publishing Group},
}

@Article{LoristNieraeth2024,
  author       = {Lorist, Emiel and Nieraeth, Zoe},
  title        = {Banach function spaces done right},
  issn         = {0019-3577},
  number       = {2},
  pages        = {247--268},
  volume       = {35},
  date         = {2024-03},
  doi          = {10.1016/j.indag.2023.11.004},
  journaltitle = {Indag. Math.},
  publisher    = {Elsevier BV},
}

@PhdThesis{Luxemburg1955,
  author      = {Luxemburg, Wilhelmus Anthonius Josephus},
  title       = {{B}anach Function Spaces},
  date        = {1955},

  institution = {Technische Hogeschool te Delft},
  url         = {https://repository.tudelft.nl/record/uuid:252868f8-d63f-42e4-934c-20956b86783f},
}

@Article{Maligranda1985,
  author       = {Lech Maligranda},
  title        = {Indices and interpolation},
  issn         = {0012-3862},
  pages        = {1-54},
  volume       = {234},
  address      = {Warszawa},
  date         = {1985},

  isbn         = {8301058188},
  journaltitle = {Diss. Math.},
  location     = {Warszawa},
  pagetotal    = {54},
  ppn_gvk      = {273704710},
  publisher    = {Inst. matematyczny, Polska Akad. nauk},
}

@Article{Mihlin1957,
  author       = {Mihlin, S. G.},
  title        = {{F}ourier Integrals and Multiple Singular Integrals},
  pages        = {143--145},
  volume       = {12},
  date         = {1957},
  journaltitle = {Vest. Leningrad Univ. Ser. Mat.},
  language     = {Russian},
  url          = {https://cir.nii.ac.jp/crid/1570009750867551360},
}

@Article{Orlicz1932,
  author       = {Orlicz, Wladyslaw},
  title        = {{\"U}ber eine gewisse {K}lasse von {R}{\"a}umen vom {T}ypus {B}},
  number       = {9},
  pages        = {207--220},
  volume       = {8},
  date         = {1932},
  journaltitle = {Bull. Int. Acad. Pol. Ser. A},
  language     = {German},
}

@Book{RaoRen1991,
  author    = {Rao, Malempati M. and Ren, Z. D.},
  publisher = {Dekker},
  title     = {Theory of {O}rlicz spaces},
  isbn      = {0824784782},
  number    = {146},
  series    = {Pure and applied Mathematics},
  date      = {1991},

  location  = {New York},
  pagetotal = {449},
  ppn_gvk   = {025373625},
}

@Article{Rolewicz1957,
  author       = {Rolewicz, Stefan},
  title        = {On a certain class of linear metric spaces},
  pages        = {471--473},
  volume       = {5},
  date         = {1957},
  journaltitle = {Bull. Acad. Polon. Sci. Cl. III},
}

@Article{Simonenko1964,
  author       = {Simonenko, I. B.},
  title        = {Interpolation and Extrapolation of Linear Operators in {O}rlicz spaces},
  number       = {4},
  pages        = {536-553},
  volume       = {105},
  date         = {1964},
  journaltitle = {Sb. Math.},
  language     = {Russian},
  url          = {https://www.mathnet.ru/eng/dan28459},
}

@InBook{TeofanovToft2022,
  author    = {Teofanov, Nenad and Toft, Joachim},
  pages     = {601--637},
  publisher = {Springer International Publishing},
  title     = {An Excursion to Multiplications and Convolutions on Modulation Spaces},
  isbn      = {9783031021046},
  booktitle = {Operator and Norm Inequalities and Related Topics},
  date      = {2022},
  doi       = {10.1007/978-3-031-02104-6_18},

  issn      = {2297-024X},
}

@Article{Toft2004a,
  author       = {Toft, Joachim},
  title        = {Continuity Properties for Modulation Spaces, with Applications to Pseudo-differential Calculus, {II}},
  issn         = {1572-9060},
  number       = {1},
  pages        = {73--106},
  volume       = {26},
  date         = {2004-08},
  doi          = {10.1023/b:agag.0000023261.94488.f4},

  journaltitle = {Ann. Glob. Anal. Geom.},
  publisher    = {Springer Science and Business Media LLC},
}

@Article{Toft2022,
  author       = {Toft, Joachim},
  title        = {Step multipliers, {F}ourier step multipliers and multiplications on quasi-{B}anach modulation spaces},
  issn         = {0022-1236},
  note         = {art. no. 109343},
  number       = {5},
  volume       = {282},
  date         = {2022-03},
  doi          = {10.1016/j.jfa.2021.109343},

  journaltitle = {J. Funct. Anal.},
  publisher    = {Elsevier BV},
}

@Article{ToftUesterNabizadehOeztop2022,
  author       = {Toft, Joachim and {\"U}ster, R{\"u}ya and Nabizadeh, Elmira and Öztop, Serap},
  title        = {Continuity Properties and {B}argmann Mappings of Quasi-{B}anach {O}rlicz Modulation Spaces},
  issn         = {1435-5337},
  number       = {5},
  pages        = {1205--1232},
  volume       = {34},
  date         = {2022-07},
  doi          = {10.1515/forum-2021-0279},

  journaltitle = {Forum Math.},
  publisher    = {Walter de Gruyter GmbH},
}
\end{document}